\documentclass[12pt]{conm-p-l}
\usepackage{}

\newtheorem{theorem}{Theorem}
\newtheorem{lemma}{Lemma}
\newtheorem{proposition}{Proposition}

\theoremstyle{definition}

\newtheorem{corollary}{Corollary}

\theoremstyle{remark}

\numberwithin{equation}{section}

\newcommand{\R}{\mathbb{R}}
\newcommand{\1}{\mathbf 1}
\newcommand{\0}{\mathbf 0}

\newcommand{\env}{\operatorname{env}}
\newcommand{\Dc}{{\mathcal D}}
\newcommand \OSup{\operatornamewithlimits \bigoplus \limits}
\newcommand \CCap{\operatornamewithlimits \bigcap \limits}

\begin{document}

\title[A Tropical Version of the Schauder
Fixed Point Theorem]{A Tropical Version of the Schauder Fixed
Point Theorem}


%    author one information
\author[G.B. Shpiz]{G.B. Shpiz}
\author[G.L. Litvinov]{G.L. Litvinov}
%%\address{Serge\u{\i} Sergeev, University of Birmingham, School of Mathematics, B15 2TT Edgbaston, Birmingham, UK}
%\curraddr{}
%%\email{sergeevs@maths.bham.ac.uk}
\thanks{This work is supported by the RFBR grant 08--01--00601 and the
joint RFBR/CNRS grant 05--01--02807}

%%    author two information
%\address{}
%\curraddr{}
%\email{}
%\thanks{}

\date{}

\begin{abstract}
A tropical version of the Schauder fixed point theorem for compact
subsets of tropical linear spaces is proved.
\end{abstract}

\maketitle

%%\newpage
\section{Introduction}

It is well-known that every continuous mapping from a compact
subset of a locally convex space to this subset has a fixed point
(the Schauder fixed point theorem \cite{1}). There exists a
correspondence (and analogy) between important, interesting, and
useful constructions and results of the traditional mathematics
over fields and analogous constructions and results over
idempotent semirings and semifields, i.e., semirings and
semifields with idempotent addition (the idempotent correspondence
principle, see \cite{2, 3, 4}). In the framework of this analogy a
tropical/idempotent version of the Schauder fixed point theorem is
proved (see Theorem 2 below).

Note that topologies in tropical/idempotent analysis do not
coincide with standard topologies in the traditional functional
analysis and collections of compact subsets of spaces of functions
do not coincide too. We shall examine applications of our results
in a separate paper. In particular, our results could be used to
prove that homogeneous (but nonlinear in general) operations in
topological idempotent linear spaces have eigenvectors. This is
closely related to asymptotic behaviour of infinite extremals in
dynamic optimization problems with infinite planning horizon, see
\cite{5}--\cite{7}. For example, it is possible to generalize the
results of \cite{5, 6} to the case of discontinuous utility
functions and kernels of Bellman operators.

\section{Basic definitions and notations}

In the present paper, we shall use some ideas and terminology from
\cite{8}--\cite{10}. We recall that an {\it idempotent semigroup}
(IS) is defined to be an additive semigroup equipped with
commutative addition $\oplus$ so that the relation $x\oplus x = x$
holds for every element $x$. Any IS can be treated as a set ordered by
the standard (partial) order: $x\preceq y$ if and only if $x\oplus
y = y$. It is easily seen that this order relation is well defined and
$x\oplus y = \sup\{ x, y\}$. For an arbitrary subset $X$ of an
idempotent semigroup, we set $\oplus X = \sup(X)$ and $\wedge X =
\inf(X)$ under the condition that the corresponding right-hand
sides exist. An {\it idempotent semiring} (ISR) is defined to be
an IS equipped with associative multiplication $\odot$ with unity
(denoted by $\1$) such that both of the distributivity relations
are satisfied. An idempotent semigroup $V$ equipped with
multiplication $\odot$ by elements from an idempotent semiring $K$
such that the relations $a\odot(b\odot x) = (a\odot b)\odot x$,
$(a \oplus b)\odot x = a\odot x \oplus b\odot x$, $a\odot(x\oplus
y) = a\odot x\oplus a\odot y$, and $\0\odot x = x\odot \0 = \0$
hold for any $a, b \in K$ and $x, y \in V$ is called an {\it
idempotent semimodule} over the idempotent semiring $K$.

The set $\R$ of all real numbers is a commutative ISR with respect
to operations $\oplus = \max$ and $\odot = +$. Denote this
semiring by $\R_{\oplus}$; we shall equip $\R_{\oplus}$ with the
standard topology of $\R$, so we shall treat $\R_{\oplus}$ as a
topological space. Note that $\R_{\oplus}$ has no zero element
$\0$ (as a semiring); if we adjoin this element then we obtain
the well-known {\it max-plus algebra} $\R_{\max} =
\R_{\oplus}\cup\{\0\} = \R_{\oplus}\cup\{-\infty\}$ or {\it
tropical algebra}. Of course, $\R_{\oplus}$ and $\R_{\max}$ have
the unity element $\1 = 0$ and the standard order in $\R_{\oplus}$
coincides with the usual one.

An idempotent semimodule over $\R_{\oplus}$ is called an {\it
idempotent $\R_{\oplus}$-space}, or $\R_{\oplus}$-{\it space}. The
semiring $\R_{\oplus}$ is an idempotent $\R_{\oplus}$-space over
itself. A homomorphism from a $\R_{\oplus}$-space $V$ to
$\R_{\oplus}$ is called a {\it linear functional} on $V$. For
arbitrary set $T$ denote by $B(T)$ the $\R_{\oplus}$-space of all
bounded mappings from $T$ to $\R_{\oplus}$ equipped with the
corresponding pointwise operations.

Let $V$ be an arbitrary partially ordered set (e.g., $V$
is an idempotent semimodule with respect to its standard
order), $a, b \in V$. We shall use the following notations for
 intervals and half-intervals:
\[
\begin{split}
[a, b] & = \{ x\in V \mid a\preceq x \preceq b\},\\
(\cdot, a] & = \{x\in V \mid x\preceq a\}, \\
[a,\cdot) & = \{x\in V \mid a\preceq x\}.
\end{split}
\]

Suppose that $X$ is a subset of an idempotent $\R_{\oplus}$-space
$V$ and $X\supset [a, b]$ for arbitrary $a \in X$ and $b\in X$.
Then we shall say that $X$ is {\it o-convex}, see \cite{11}. We
shall say that a topology on $V$ is {\it locally o-convex} if
every element $x\in V$ has a basis of o-convex neighborhoods.
Suppose that $V$ is  equipped with an o-convex topology such that
for arbitrary $v\in V$ the mapping $r \mapsto r\odot v$ (from
$\R_{\oplus}$ to $V$) is continuous and the half-intervals
$(\cdot, v]$ and $[v, \cdot)$ are closed. Then we shall say that
$V$ is a {\it topological $\R_{\oplus}$-space}.

Denote by $V^*$ the set of all continuous linear (over
$\R_{\oplus}$) functionals on $V$. The set $V^*$ is an
$\R_{\oplus}$-space with respect to the corresponding pointwise
operations. We shall say that $V$ is {\it regular}, if for every
$x, y \in V$, $x\neq y$, there exists a functional $w\in V^*$ such
that $w(x) \neq w(y)$. The topology generated by the basis of all
sets of the form $\{x\in V\mid a< w(x) < b\}$ for $a, b\in
\R_{\oplus}$, $w\in V^*$ will be called a $\oplus$-{\it weak
topology}.

Suppose that $V$ and $W$ are topological $\R_{\oplus}$-spaces and
$f$ is a mapping from $V$ to $W$. This mapping is called
$\oplus$-{\it weakly continuous} if for every $w\in W^*$ the
mapping $wf : V\to \R_{\oplus}$ is continuous.

Suppose that $V$ is an $\R_{\oplus}$-space and $x, y \in V$. We
shall write $x\gg y$ if there exists an element $r >\1$ ($r\in
\R_{\oplus}$) such that $r\odot y \preceq x$. Define subsets
$\Dc_x(r)\subset V$ by the following formula:
\[
\begin{split}
\Dc_x(r) & = \{ y\in V\mid r\odot x\gg y \gg r^{-1}\odot x\}\\
& = \{ y\in V\mid r\odot x\gg y \mbox{ and } r\odot y\gg x\}.
\end{split}
\]

The topology generated by the basis of all sets of the form
$\Dc_x(r)$ for $r > \1$ will be called {\it uniform}. The uniform
topology is metrizable. The corresponding metric can be defined,
e.g., by the formula:
\[
d(x, y) = \arctan(\inf\{r\in \R_{\oplus}\mid r^{-1}\odot x\preceq
y\preceq r\odot x\}).
\]
In the space $B(X)$ of all bounded real functions defined on a set
$X\neq\emptyset$ the uniform topology is defined by the metric
\[
d(f, g) = \sup_{x\in X}\limits\mid f(x) - g(x)\mid.
\]

\section{Topological $\R_{\oplus}$-spaces}

\begin{lemma}
Let $V$ be an arbitrary $\R_{\oplus}$-space. Then the sets of the
form $\Dc_v(l)$ for $l > \1$ form a basis of neighborhoods of the
point $v\in V$ with respect to the uniform topology.
\end{lemma}

\begin{proof}
It is necessary to check that for every $y\in \Dc_x(l)$ there
exists $r > \1$ such that $\Dc_y(r)\subset \Dc_x(l)$. We have
$l\odot x\gg y$ and $l\odot y\gg x$, so there exists $p > \1$ such
that $l\odot x\succeq p\odot y$ and $l\odot y\succeq p\odot x$. We claim that
any $r$ such that $\1 < r < p$ is good enough.
For any $z\in\Dc_y(r)$, we have $z\gg r^{-1}\odot y$ and $r\odot y\gg z$. 
It follows that
$l\odot z\gg l\odot r^{-1}\odot y\gg
p^{-1}\odot l\odot y\succeq x$, and on the other hand $l\odot
x\succeq p\odot y\gg r\odot y\gg z$. Thus $z\in \Dc_x(l)$. So
$\Dc_y(r)\subset \Dc_x(l)$ because $z\in \Dc_y(r)$ is an arbitrary
element.
\end{proof}

\begin{proposition}
Let $V$ and $W$ be $\R_{\oplus}$-spaces. Suppose that $f : V\to W$
is a nondecreasing mapping such that $f(r\odot v)\preceq r\odot
f(v)$ for all $r\succeq \1$, $v\in V$. Then the mapping $f$ is
continuous with respect to the uniform topology.
\end{proposition}

\begin{proof}
From Lemma 1 it follows that it is sufficient to prove that the
preimage of $\Dc_{f(x)}(l)$ contains a neighborhood of $x$. In
fact the preimage of $\Dc_{f(x)}(l)$ contains $\Dc_x(l)$. Indeed,
suppose that $y\in \Dc_x(l)$. Then $l\odot x\succeq r\odot y$ and
$l\odot y\succeq r\odot x$ for some $r\in (\1, l)$. So we have
$l\odot r^{-1}\odot f(x)\succeq f(l\odot r^{-1}\odot x)\succeq
f(y)$ and $l\odot r^{-1}\odot f(y)\succeq f(l\odot r^{-1}\odot
y)\succeq f(x)$. Therefore, we have $l\odot f(x)\gg f(y)$ and
$l\odot f(y)\gg f(x)$, that is $f(y)\in \Dc_{f(x)}(l)$ and the
proposition is proved.
\end{proof}

\begin{corollary}
For every $\R_{\oplus}$-space, both addition $\oplus$ and 
multiplication by a number are continuous with respect to the uniform topology.
\end{corollary}

\begin{proposition}
For every topological $\R_{\oplus}$-space $V$ the following
statements hold:
\begin{enumerate}
\item The topology of $V$ is majorized by the uniform topology of
$V$, i.e. every open subset of $V$ is open for the uniform
topology.
\item The space $V$ is a topological $\R_{\oplus}$-space with
respect to the uniform topology.
\end{enumerate}
\end{proposition}

\begin{proof}
\quad
\begin{enumerate}
\item Suppose that $x\in V$ and $U$ is a neighborhood of $x$. Let us
show that $U$ contains a neighborhood of $x$ with respect to the
uniform topology. The space $V$ is locally o-convex, so we can
assume that $U$ is an o-convex set. The mapping $r\mapsto r\odot
x$ is a continuous mapping from $\R_{\oplus}$ to $V$, so $U$
contains $r\odot x$ and $r^{-1}\odot x$ for some $r > \1$. Thus
$\Dc_x(r) \subset U$ and the statement is proved.
\item It is obvious, that the uniform topology is o-locally convex
and mapping $r\mapsto r\odot x$ is continuous. From the statement
(1) it follows that all the half-intervals of the form $(\cdot,
a]$ and $[a, \cdot)$ are closed. Thus the proposition is proved.
\end{enumerate}
\end{proof}

The proof of the following proposition is straightforward.

\begin{proposition} If $T$ is a finite nonempty set, the uniform
topology for $B(T)$ is the topology of pointwise convergence, that
is the usual topology of the Euclidean space.
\end{proposition}

Below we shall equip $B(T)$ with the uniform topology.

\begin{proposition}
Suppose that $V$ is a regular topological $\R_{\oplus}$-space.
Then $V$ is a topological $\R_{\oplus}$-space with respect to its
$\oplus$-weak topology.
\end{proposition}

\begin{proof}
It is easy to see that the $\oplus$-weak topology is o-convex and
the mapping $r\mapsto r\odot a$ is continuous under this topology.
Let us show that the sets of the form $(\cdot, a]$ and $[a,
\cdot)$ are closed. Suppose that $x$ is an element of the closure
of $(\cdot, a]$ relative to the $\oplus$-weak topology. Then
$w(x)\preceq w(a)$ for every functional $w\in V^*$, so $w(a\oplus
x) = w(a)$. Hence $a\oplus x = a$ because $V$ is regular, so
$x\preceq a$, that is $x\in (\cdot, a]$ and the set $(\cdot, a]$
is closed. For the set $[a, \cdot)$ the proof is similar. So the
proposition is proved.
\end{proof}

\section{Convex subsets in topological $\R_{\oplus}$-spaces}

Let $V$ be an idempotent semimodule over an idempotent semiring
$K$, $X$ a subset of $V$, and $p : X \to K$ a function such that
$\oplus p(x)=\1$. The element $\OSup _{x \in X} (p(x)\odot x)$ is
called a $\oplus$-{\it convex combination} of all the elements
$x\in X$. A subset $X\subset V$ is called $\oplus$-{\it convex} if
$X$ contains every $\oplus$-convex combination of elements of
every finite subsets of $X$, see \cite{12}--\cite{15} for general
definitions and constructions of this type. We shall say that a
subset $X\subset V$ is {\it $a$-convex} if there exists the
$\oplus$-convex combination  $\OSup _{x \in X} (p(x)\odot x)\in X$
for each function  $p: X\to K$ such that  $\oplus p(X)=\1$.

From these definitions it follows that every $\oplus$-convex set
is a subsemigroup with respect to the idempotent addition $\oplus$ and
every $a$-convex set is bounded with respect to the standard order
(see Section 2 above). Of course, every $a$-convex set is
$\oplus$-convex.

\begin{proposition}
Let $V$ be a topological $\R_{\oplus}$-space, $X$ its compact
subsemigroup. Then there exists $\oplus X\in X$.
\end{proposition}

\begin{proof}
For  $v\in X$ we set  $X(v)=\{x\in X\mid v\preceq x\} =
[v,\cdot)\cap X$. For each finite subset $A\subset X$ we have
$\oplus A\in \CCap_{v\in A}X(v)$, so the collection $\{X(v)\}$ is
a centered family of closed subsets of $X$. The set $X$ is
compact, so there exists an element $x\in \CCap_v X(v)$. By our
construction $x\succeq X$ and $x\in X$, so $x = \oplus X$. Thus
there exists $\oplus X$ and $\oplus X\in X$. The proposition is
proved.
\end{proof}

\begin{corollary}
Under the conditions of Proposition 5 for arbitrary subset $Y$ of
$X$ there exists the sum  $\oplus Y$ and this sum is an element of
the intersection of all closed subsubsemigroups containing  $Y$.
Moreover, for each $w\in V^*$ we have $w(\oplus Y) = \oplus w(Y)$.
\end{corollary}

\begin{proof}
Denote by $\widehat Y$ the intersection of all closed
subsemigroups containing  $Y$; the set $\widehat Y$ is a compact
semigroup, so it is possible to apply Proposition 5 and the first
statement of the corollary is proved. To prove the second
statement it is sufficient to prove that $\oplus Y=\oplus
(\widehat Y)$. But it follows from the obvious fact: the
set $\{x\in V\mid x\preceq b\}$ is a closed subsemigroup for each
$b\in V$.
\end{proof}

\begin{proposition}
Let $V$ be a topological $\R_{\oplus}$-space. Each compact
$\oplus$-convex subset $X$ of $V$ is a-convex.
\end{proposition}

\begin{proof}
From Proposition 5 it follows that for each subset $Y$ of $X$ the
sum $\oplus Y$ exists and $\oplus Y\in X$. Let  $p: X\to K$ is a
function such that  $\oplus p(X)=\1$. Denote by $\widehat p$ the
convex combination $\OSup_{x\in X}(p(x)\odot x)$. If  $p(v)=\1$
for an element $v\in X$, then  $\widehat p = \OSup_{x\in
X}(p(x)\odot x\oplus v)$. For each $x\in X$ the element $p(x)\odot
x\oplus v$ belongs to $X$, hence $\widehat p\in X$.  If the
function $p$ does not reach its maximum, then for an arbitrary
number $r < \1$ we set $p_r(x) = r^{-1}\odot (p(x)\wedge r)$  and
 $\widehat {p_r} = \OSup_{x\in X}(p_r(x)\odot x)$; recall that
 $p(x)\wedge r = \inf\{ p(x), r\} = \min\{p(x), r\}$. By
 construction we have $\oplus p_r(X) = \OSup_{x\in X}p_r(x) = \1$
 and $ p_r(v)=\1$ for an element  $v\in X$. By construction we have
  $\widehat p\preceq \widehat {p_r}\preceq r^{-1}\odot\widehat p$, so
  $\widehat p_r$ converges to $\widehat p$ with respect to the
  uniform topology as $r$ tends to $\1$. From Proposition 2 it follows that
  $\widehat p_r$ converges to  $\widehat p$ in $V$ as $r$ tends to
  $\1$, so $\widehat p\in X$ because $X$ is compact. The
  proposition is proved.
\end{proof}

Let $V$ be a topological $\R_{\oplus}$-space, $X$ its subset. We
set
\[
\env(X) = \{y\in  V \mid (\exists x\in X, r\in \R_{\oplus})\;
r\odot x\preceq y\}.
\]
It is clear that $\R_{\oplus}\odot X\subset\env(X)$.

\begin{proposition}
Let $V$ be a topological $\R_{\oplus}$-space, $X$ its a-convex
subset. Then there exists a mapping $\pi : \env(X)\to X$ such that
$\pi$ is continuous with respect to the uniform topology, $\pi(x)
= x$, and $\pi (r\odot x) \in \R_{\oplus}\odot x$ for all $x\in
X$, $r\preceq\1$.
\end{proposition}

\begin{proof}
Set
\[
M = \{(x,y)\in X\times \env(X) \mid (\exists\; r\preceq \1)\;
r\odot x\preceq y\}.
\]
For $(x, y)\in M$ we set
\[
\begin{split}
r_x(y) & = \sup\{ r\preceq\1 \mid r\odot x\preceq y\},\\
m(y) &=  \OSup_{(z,y)\in M}r_z (y),\\
n_x(y) &= m(y)^{-1}\odot r_x (y),\\
p(y) &= \OSup_{(z,y)\in M}r_z(y)\odot z,\\
\pi (y) &= \OSup_{(z,y)\in M}n_z(y)\odot z = m(y)^{-1}\odot p(y).
\end{split}
\]

By construction we have $m(y)\preceq \1$ and
 $\OSup_{ x\in X}n_x(y) =\1$ for $y\in \env(X)$.
 Since $\pi(y)$ is a $\oplus$-convex combination of elements of $X$
 the element $\pi(y)$ belongs to $X$. We have $m(y)\odot n_x(y)
 \odot x = r_x(y)\odot x\preceq y$; hence
  $m(y)\odot \pi (y)\preceq y$. By  construction we have
$m(x) = r_x(x) = n_x(x)=\1$ for $x\in X$, so $x = n_x(x)\odot
x\preceq \pi(x)\preceq x$. Thus $\pi(x) = x$ and $\pi$ is a
retraction $\env(X) \to X$.

Let us prove that $\pi$ is continuous. Since $\pi(x) =
m(x)^{-1}\odot p(x)$ and the multiplication by coefficients is a
continuous mapping $\R_{\oplus}\times V\to V$ with respect to the
uniform topology it is sufficient to show that  $m$ and $p$ are
continuous for the uniform topology. This statement follows from
Proposition 1. Thus Proposition 7 is proved.
\end{proof}

Note that similar constructions were used in \cite{12, 14}.

\begin{proposition}
Suppose that $T$ is a finite nonempty set, $X$ is a compact
$\oplus$-convex subset of $B(T)$ and $Y$ is the traditional
(usual) convex hull of $X$ in the Euclidean space $B(T)$. Then
there exists a continuous mapping $\pi : Y\to X$ such that $\pi(X)
= X$.
\end{proposition}

\begin{proof}
From Proposition 6 it follows that the set $X$ is $a$-convex. It is
clear that $Y\subset \env(X) = B(T)$, so it is possible to apply
Proposition~7.
\end{proof}

The following theorem is a tropical/idempotent version of the
Brauer fixed point theorem.

\begin{theorem}
Suppose that $T$ is a nonempty finite set, $X$ is a compact
$\oplus$-convex subset of $B(T)$, and $f$ is a continuous mapping
from $X$ to $X$. Then $f$ has a fixed point.
\end{theorem}

\begin{proof}
Suppose that $Y$ is the traditional convex hull of $X$ in the
Euclidean space $B(T)$ and $\pi : Y\to X$ is the continuous
mapping (retraction) discussed in Proposition 7. Then $Y$ is a
compact convex subset (in traditional sense) in the Euclidean
space $B(T)$. From the Brauer fixed point theorem it follows that
the mapping $f\pi$ has a fixed point $x\in Y$. In fact $x\in X$
because $f\pi(Y) \subset X$. Hence $\pi(x) = x$ and $x$ is a fixed
point for the mapping $f$. Theorem 1 is proved.
\end{proof}

\section{Main results}

Suppose that $V$ is a topological $\R_{\oplus}$-space, $X$ is a
compact $\oplus$-convex subset of $V$, $T$ is a nonempty finite
subset of $V^*$. Consider the mapping
\[
i : V\to B(T)
\]
defined by the formula $i(v) : v\mapsto t(v)$, where $v\in V$ and $t\in T$.
Obviously, the mapping $i$ is linear and continuous.

\begin{lemma}
There exists a continuous mapping $p : B(T)\to X$ such that
$i(p(f)) = f$ for each $f\in i(X)$.
\end{lemma}

\begin{proof}
For $f \in B(T)$ we set $p(f) = \OSup_{x\in X, i(x)\preceq f}x$.
By Proposition 6, $X$ is $a$-convex. Hence $p(f)\in X$, 
and for $r\succeq \1$ and $f\in B(T)$ we have
 $r^{-1}p(r\odot f)\oplus p(f)\in X$.
We have $i(p(f))\preceq f$.
Using this inequality and the linearity of $i$,
we obtain $i(r^{-1}p(r\odot f)\oplus p(f))
 \preceq r^{-1}\odot r\odot f\oplus f = f$. So we have
$r^{-1}p(r\odot f)\preceq p(f)$, that is $p(r\odot f)\preceq
r\odot p(f)$. By construction the mapping $p$ is nondecreasing; so
from Proposition 1 it follows that $p$ is continuous for the
uniform topology on $V$. From Proposition 2 it follows that $p$ is
continuous with respect to the initial topology (on $V$) which is
weaker than the uniform topology. Thus Lemma 2 is proved.
\end{proof}

\begin{proposition}
Suppose that $V$ is a topological $\R_{\oplus}$-space, $X$ is its
compact $\oplus$-convex subset, $f$ is a $\oplus$-weakly
continuous mapping from $X$ to $X$, and $T$ is a nonempty finite
subset of $V^*$. Then there exists an element $x\in X$ such that
$w(x) = w(f(x))$ for each $w\in T$.
\end{proposition}

\begin{proof}
Suppose that $i : V\to B(T)$ is the mapping defined in the
beginning of this section (before Lemma 2) and $U = i(X)$. 
As $i$ is linear and continuous, $U$ is a compact
$\oplus$-convex subset of $B(T)$. From Lemma 2 it follows that
there exists a continuous mapping $p :U\to X$ such that $i(p(w)) =
w$ for every functional $w\in U$. The formula $g(w) = i(f(p(w)))$,
where $w \in U$, generates a mapping $g : U\to U$. This mapping is
continuous, $U$ is a compact $\oplus$-convex subset of $B(T)$. So,
by Theorem 1, the mapping $g$ has a fixed point $u\in U$. Set $x =
p(u)$. Since $i(f(p(u))) = u$, we have $i(x) = i(p(u)) = u =
i(f(p(u))) = i(f(x))$, that is $w(x) = w(f(x))$ for each $w\in T$.
The proposition is proved.
\end{proof}

The following theorem is a tropical/idempotent version of the
Schauder fixed point theorem.

\begin{theorem}
Suppose that $V$ is a regular topological $\R_{\oplus}$-space, $X$ is its
compact $\oplus$-convex subset, and $f$ is a $\oplus$-weakly
continuous mapping from $X$ to $X$. Then $f$ has a fixed point
$x\in X$.
\end{theorem}

\begin{proof}
For every nonempty finite subset $T$ of $V^*$ we define a set
$S(T)$ by the formula
\[
S(T) = \{ x\in X \mid w(x) = w(f(x)) \mbox{ for each } w\in T\}.
\]
 By construction this set is closed; by Proposition 9 it is
nonempty. Obviously, $S(T_1)\cap S(T_2) = S(T_1\cup T_2)$; so the
family of all sets of the form $S(T)$ is a centered family of
closed subsets of the compact set $X$. Hence this family has a
nonempty intersection. Let $x$ be an element of this intersection.
By construction we have $w(x) = w(f(x))$ for all $w\in V^*$.
Therefore, $f(x) = x$ because $V$ is regular. Theorem 2 is proved.
\end{proof}

\section{Acknowledgement}
The authors are grateful to S.~N.~Sergeev for 
his kind help, and to the anonimous referee
for careful reading and useful remarks on the previous
version of the paper.

\end{document}